\documentclass[a4paper]{my_elsarticle}

\usepackage{amsmath,amssymb,amsfonts,amsthm}

\newtheorem{example}{Example}
\newtheorem{remark}[example]{Remark}

\begin{document}

\begin{frontmatter}

\title{Exponentials and Laplace transforms\\ on nonuniform time scales}

\author{Manuel Ortigueira}
\ead{mdo@fct.unl.pt}
\address{CTS-UNINOVA, Department of Electrical Engineering,
Faculty of Science and Technology,\\
Universidade Nova de Lisboa, Portugal}

\author{Delfim F. M. Torres}
\ead{delfim@ua.pt}
\address{Center for Research and Development in Mathematics and Applications (CIDMA),\\
Department of Mathematics, University of Aveiro, 3810--193 Aveiro, Portugal}

\author{Juan Trujillo}
\ead{jtrujill@ullmat.es}
\address{Universidad de La Laguna,
Departamento de An\'{a}lisis Matem\'{a}tico,\\
38271 La Laguna, Tenerife, Spain}


\begin{abstract}
We formulate a coherent approach to signals and systems theory on time
scales. The two derivatives from the time-scale calculus are used, i.e.,
nabla (forward) and delta (backward), and the corresponding eigenfunctions, the
so-called nabla and delta exponentials, computed. With these exponentials, two
generalised discrete-time Laplace transforms are deduced and their properties
studied. These transforms are compatible with the standard Laplace and $Z$
transforms. They are used to study discrete-time linear systems defined
by difference equations. These equations mimic the usual continuous-time
equations that are uniformly approximated when the sampling interval becomes small.
Impulse response and transfer function notions are introduced. This implies
a unified mathematical framework that allows us to approximate the classic
continuous-time case when the sampling rate is high or to obtain the standard
discrete-time case, based on difference equations, when the time grid becomes uniform.
\end{abstract}

\begin{keyword}
time-scale calculus 
\sep exponentials
\sep generalized Laplace and $Z$ transforms
\sep systems theory
\sep fractional derivatives.

\MSC[2010] 26E70 \sep 44A10 \sep 65T50.

\end{keyword}
\end{frontmatter}


\section{Introduction}

The analysis of nonuniformly sampled data is a very important task having large spread
application in fields like astronomy, seismology, paleoclimatology, genetics and laser 
Doppler velocimetry \cite{BaSt}. The ``jitter'' in Telecommunications is a well known 
problem \cite{GlGr}. Very interesting is the heart rate variability of the signal 
obtained from the ``R'' points \cite{Pe}.
Traditionally, most interesting techniques for dealing with this kind of signals 
pass by interpolation, to obtain a continuous-time signal that is analysed 
by current methods \cite{AlGr,FeGr:1,FeGr:2}. An alternative approach proposed 
in \cite{Gr} allows a conversion from irregular to regular samples, 
maintaining the discrete-time character. Other approaches, include the study 
of difference equations with fractional delays \cite{Bal14,Orti00}. 
However, no specific tools for dealing directly with such signals were developed. 
In particular, no equivalent to the Laplace or $Z$ transforms were proposed. 
As well-known, the use of Laplace and $Z$ transforms,
to solve differential and difference linear equations,
is very common in almost all scientific activities \cite{PM}.
Normally, uniform time scales are used, but frequent applications
use nonuniform scales. This makes important to obtain generalisations
of such transforms for other kinds of scales. Some attempts have been made
\cite{A09,AH90,AH90b,BG10,BG10b}, but let us an unsatisfactory feeling:
they are not true generalisations of the classic formulations. The main
difficulty is in the starting point, i.e., the exponentials
used to define the transforms. Usually, on time scales, causal exponentials
are used instead of two-sided exponentials \cite{C12,HAQ10,HAQ12,MT09}.
On the other hand, no correct interplay between nabla and delta derivatives
and exponentials and transforms has been established. Such interplay 
was stated for fractional derivatives in the recent paper \cite{OrCoTr15}. 

Here, in a first step, we clarify nabla and delta definitions and their meaning 
and relation with causality. With nabla (causal) and delta (anti-causal) 
derivatives, we define corresponding linear systems. Each concept is used 
to define two exponentials over the whole time scale and not only above 
or below a given time reference. With each exponential, a given transform 
is defined. We start from the inverse transform and only later we define 
the direct transform. With the nabla exponential, we define the inverse 
nabla transform through a Mellin-like integral on the complex plane. 
The direct nabla transform is defined with the help of the delta exponential. 
For the delta transform we reverse the exponentials.
Having defined the exponentials, we study the question of existence,
arriving to the concept of region of convergence. The unicity of the 
transforms is also investigated. This lead us to generalise the convolution
and correlation concepts with the help of equivalent time scales.
The concept of transfer function, as the eigenvalue corresponding
to the respective exponential, is introduced. The inverse Laplace transform
is the so called impulse response, i.e., the response
of the system when the input is a delta function.
We consider also the conversion from one time scale to another one that
is equivalent to it. In passing, we prove the existence
of no periodicity and, consequently, the inability to define Fourier
transforms and series. This is the main drawback of the theory.


\section{On the calculus on time scales}

A powerful approach into the continuous/discrete unification/generalization 
was introduced by Aulbach and Hilger through the calculus 
on \textit{measure chains} \cite{AH90,H90}. However, the main popularity 
was gained by the calculus on \textit{time scales} \cite{AH90b,BP01,BP03}. 
These are nonempty closed subsets $\mathbb{T}$ of the set $\mathbb{R}$
of real numbers, particular cases of measure chains. We remark that the name may
be misleading, since the term \emph{scale} is used in Signal Processing
with a different meaning. On the other hand, in many problems 
we are not dealing with time.

Let $t$ be the current instant. Using the language of the time-scale calculus,
the previous next instant is denoted by $\rho(t)$. Similarly, the next following 
point on the time scale $\mathbb{T}$ is denoted by $\sigma(t)$. One has
\begin{equation*}
\rho(t) = t - \nu(t),
\quad  \sigma(t) = t + \mu(t),
\end{equation*}
where $\nu(t)$ and $\mu(t)$ are called the graininess functions.
We avoid here the use of the terms forward and backward graininess,
since they have meanings that are different from the ones used
in Signal Processing applications. Let us define
\begin{equation*}
\nu^0(t) := 0, \quad
\nu^{n}(t) := \nu^{n-1}(t) + \nu\left(t-\nu^{n-1}(t)\right),
\quad n \in \mathbb{N},
\end{equation*}
and $\rho^0(t) := t$, $\rho^n(t) := \rho\left(\rho^{n-1}(t)\right)$,
$n \in \mathbb{N}$. Note that $\nu^1(t) = \nu(t)$ and $\rho^1(t) = \rho(t)$.
When moving into the past, we have
\begin{equation*}
\begin{split}
\rho^0(t) &= t = t - \nu^0(t),\\
\rho^1(t) &= \rho(t) = t-\nu(t) = t - \nu^1(t),\\
\rho^2(t) &= \rho(\rho(t)) = \rho(t)-\nu(\rho(t))
= t-\nu(t)-\nu(t-\nu(t)) = t - \nu^2(t),\\
&\vdots\\
\rho^n(t) &= t - \nu^n(t).
\end{split}
\end{equation*}
Moving into the future, the definitions and results are similar:
\begin{gather*}
\mu^0(t) := 0, \quad \mu^{n}(t) := \mu^{n-1}(t) + \mu\left(t+\mu^{n-1}(t)\right),\\
\sigma^0(t) := t, \quad \sigma^n(t) := \sigma\left(\sigma^{n-1}(t)\right),
\end{gather*}
$n \in \mathbb{N}$, and we have $\sigma^n(t) = t + \mu^n(t)$.

\begin{example}
Let $\mathbb{T} = h \mathbb{Z}$, $h > 0$.
In this case one has $\mu^n(t) = \nu^n(t) = n h$,
$\sigma^n(t) = t + n h$, and $\rho^n(t) = t - n h$, $n= 0, 1, 2, \ldots$
\end{example}

With our notations, a time scale of isolated points is written in the form
$$
\mathbb{T} = \{\ldots, \rho^{m}(t), \ldots, \rho^2(t), \rho^1(t), t,
\sigma^1(t), \sigma^2(t), \ldots \sigma^n(t), \ldots\}.
$$
This means that we can refer to all instants in the time scale $\mathbb{T}$
with respect to only one point $t$, taken as the reference. For example,
if $\mathbb{T} = \mathbb{Z}$, such instant is usually taken as $t=0$.
We must be careful with the meaning of differences like $t - \tau$, because
they can make no sense. Assume that $t$ and $\tau$ are in $\mathbb{T}$ with $t > \tau$.
This means that there exits an integer $N$ such that $\tau = \rho^N(t) = t - \nu^N(t)$
and $t = \sigma^N(\tau) = \tau + \mu^N(\tau)$. It follows that
$t - \tau = \nu^N(t) = \mu^N(\tau)$. Let us now define
$$
\nu_n(t) := \rho^{n-1}(t) - \rho^{n}(t),
\quad
\mu_n(t) := \sigma^{n}(t) - \sigma^{n-1}(t),
$$
as the $n$th values of the graininess functions from $t$.
One can then see that the difference $t - \tau $ above
is equal to the sum of the graininess
values $\mu_n$ to go from $\tau$ to $t$,
$$
t - \tau = \mu^N(\tau) = \mu(\tau)+\mu(\tau+\mu(\tau)) + \cdots
= \sum\limits_{n = 1}^{N} \mu_n(\tau),
$$
or, which is the same, the sum of the graininess values
$\nu_n$ to go from $t$ to $\tau$,
$$
t - \tau = \nu^N(t) = \nu(t)+\nu(t-\nu(t)) + \cdots
= \sum\limits_{n = 1}^{N} \nu_n(t).
$$

In this work we consider time scales $\mathbb{T}$ defined
by a set of discrete instants $t_n$, $n \in \mathbb{Z}$,
and by the corresponding graininess functions.
These instants are isolated points or the consecutive boundary points
defining a closed interval, in which the graininess functions are null.
In this case, we define the \emph{graininess interval} to be the width
of the interval. Let
$\mathbb{T} =  \bigcup_{n=1}^{\infty} \left[ 2n, 2n+1 \right]$.
Accordingly to what we just said, the graininess functions are zero
inside the intervals starting at the even integers $2n$
and finishing at the odd $2n + 1$. Moreover, $\mu(2n + 1) = 1$.

To fix ideas and have a useful framework, we assume that our domain
is a set formed by joining two sets: a set of isolated points and a set
made by closed intervals. Now we consider that $\mathbb{T}$
is the set of the given discrete set and the extreme points of the intervals.
We consider the discrete set $\mathbb{T}=\left\lbrace t_n: \;
n\in \text{$\mathbb{Z}$} \right\rbrace$ and attach to each point a label with:
isolated, left dense, or right dense. Each left dense point (say $t_k$)
represents all the points $t_{k-1}+h$ with the convention of inserting
a limit computation when computing the nabla derivative and transforming
a summation into an integral. If the point is right dense the situation
is similar for the interval defined by $t_{k}-h$ and for the delta derivative.
We define a direct graininess
\[
t_n = t_{n-1} + \nu_n, \quad n \in \mathbb{Z},
\]
and reverse graininess
\[
t_n = t_{n+1} - \mu_n, \quad n \in \mathbb{Z},
\]
where we avoid representing the reference instant $t_0$.
These definitions are suitable for dealing with some of the most interesting time
scales we find in practice. However, we have some difficulties in doing some kind
of manipulations that are also very common. Let us consider a time scale defined
on $\mathbb{R}$ and unbounded when $t \rightarrow \pm \infty$.
For example, a time scale defined by
\[
t_n = nT + \tau_n,
\quad n \in \mathbb{Z},
\quad T>0,
\quad \left| \tau_n \right| < \frac{T}{2},
\]
which we can call \textquotedblleft almost
linear sequence\textquotedblright \ \cite{Orti01}, or
\[
t_n = t_{n-1} + \tau_n,
\quad n \in \mathbb{Z},
\quad \tau_n > 0.
\]
We can consider intervals with null graininess. These time scales lead to the
generalized Laplace and $Z$ transforms. Now, consider a power transformation
of the time scale:
\[
\theta = q^t.
\]
It is clear that the new time scale $\Theta$ may have zones with null
graininess that cannot be treated as referred above. In doing limit computations,
when the graininess is null, in the above case we substitute it by an $h$
that decreases to zero. In this case we must substitute by $t(1-q) $ or
$t(q^{-1}-1)$, with $q < 1$. Besides, the new time scale is in $\mathbb{R}^+$.
We will have time scales like
\[
\theta_n = \theta^n,
\quad n \in \mathbb{Z}, \quad \theta>0,
\]
or
\[
\theta_n = \theta_{n-1} \cdot \tau_n,
\quad n \in \mathbb{Z},
\quad \tau_n > 0.
\]
Now, consider a new time scale $\overline{\mathbb{T}}$,
called a \emph{super time scale}, resulting from $\mathbb{T}$ through
\[
\overline{\mathbb{T}} = \left\lbrace \bar{t} \in \mathbb{R} :
\bar{t} = t_n - t_k, \: t_n, t_k\in \mathbb{T} \right\rbrace.
\]
If the origin $0$ is in $\mathbb{T}$, then $\overline{\mathbb{T}}$
contains $\mathbb{T}$. If not, join both and continue to call it
$\overline{\mathbb{T}}$ that will surely have $0$ as an element.
Let us consider all possible sub time scales that we construct
from $\overline{\mathbb{T}}$ and that we can put in bijection with $\mathbb{T}$.
Call them \emph{equivalent grids}; later we will give some justification for this name.


\section{On the Laplace and $Z$ transforms}

The importance of Laplace and $Z$ transforms in the study of linear
time invariant systems is unquestionable. Many books enhance such fact.
However, most of them do not show clearly why they are important and this led
to misinterpretations that conditioned their generalizations to time scales
\cite{BG07,MGD08}.


\subsection{Relations between transforms and linear time invariant systems}

Let us first consider the continuous time case.
In this case, the input-output relation for linear
time invariant (LTI) systems is given by convolution \cite{Rob03}:
\begin{equation}
\label{LTZ_2}
y(t) = \int\limits_{-\infty}^{+\infty} h(\tau) x(t - \tau) \mathrm{d}\tau,
\end{equation}
where $x(t) $ is the input signal and $h(t)$ is the impulse response.
Let $x(t) = e^{st}$, $t \in \mathbb{R} $. As it is easy to see, the output
is given by $y(t) = H(s)e^{st}$ with $H(s)$ given by
\begin{equation*}
H(s) = \int\limits_{-\infty}^{+\infty} h(t) e^{-st} \mathrm{d}t,
\end{equation*}
which is called the \emph{transfer function}. As we see, $H(s)$ is the bilateral
(two-sided) Laplace transform of the impulse response of the system. We showed
that the exponential, defined in $\mathbb{R}$, is the eigenfunction
of a continuous-time LTI system. This shows why we must use
the two-sided Laplace transform, although the one-sided has greater popularity.
However, the exponential used in this transform is not an eigenfunction
of a LTI system. It is currently used because of the initial conditions,
but it fails in the fractional case \cite{Orti11}.
Moreover, the two-sided has as special case the Fourier transform:
we only have to make $s=i\omega$. In the following,
we always consider the two-sided Laplace transform.\footnote{We easily obtain
the one-sided Laplace transform by multiplying the function
to be transformed by the Heaviside unit step.}
Now we consider the Laplace transform of a signal resulting from the uniform
sampling of a given function, $f(t)$, defined in $\mathbb{R}$. Let the ideal
sampler be the comb defined by a sequence of Dirac deltas \cite{Orti01}:
\begin{equation*}
p(t)= \sum\limits_{n=-\infty}^{+\infty} \delta(t - nh),
\end{equation*}
where $h > 0$ is the sampling interval. The ideal sampled function is,
by the properties of the Dirac delta,
\begin{equation*}
f_p(t) = \sum\limits_{n=-\infty}^{+\infty} f(nh)\delta(t - nh).
\end{equation*}
Its Laplace transform is given by
\begin{equation*}
F_p(s) = \sum\limits_{n=-\infty}^{+\infty} f(nh) e^{-shn},
\end{equation*}
which can be considered the discrete Laplace transform of function
$f_n = f(nh)$. With a change of variable $z = e^{sh}$, we obtain the $Z$ transform
\begin{equation}
\label{LTZ_7}
F(z) = \sum\limits_{n=-\infty}^{+\infty} f_n z^{-n}.
\end{equation}
As it is easy to verify, the discrete exponential $z^n$ is also the eigenfunction
of a discrete-time linear system defined by linear difference equations.
It is important to remark that $z^{-1}$ represents a delay in time. Frequently,
we convert continuous-time systems to discrete-time systems by doing a conversion from
$s$ to $z$ through the Euler transformation
\begin{equation*}
s = \frac{1 - z^{-1}}{h}
\end{equation*}
or the bilinear transformation
\begin{equation*}
s = \frac{2}{h}\frac{1 - z^{-1}}{1 + z^{-1}}.
\end{equation*}
The former seems more natural, due to its relation to the incremental ratio
used to define the derivative,
\begin{equation}
\label{LTZ_9}
f'(t) = \lim\limits_{h\rightarrow 0}\frac{f(t) - f(t-h)}{h},
\end{equation}
and to the approximation $z^{-1} \approx 1 - sh$ for small $h$.
However, the bilinear transformation is considered to be more useful,
since on converting the left $s$ half plane into the unit circle in the $z$
plane allows the design of discrete-time systems from continuous-time prototypes.
From these considerations, we must restate that:
\begin{itemize}
\item The two-sided transforms appear naturally as consequence
of the fact that exponentials are eigenfunctions of LTI systems.
\item The variable $s$ in the Laplace transform represents the derivative,
while $z^{-1}$ in the $Z$ transform represents a delay and \eqref{LTZ_7}
establishes a relation between them.
\item The derivative in \eqref{LTZ_9} is preferable because it uses
the present and past values in contra-position to
\begin{equation}
\label{LTZ_10}
f'(t) = \lim\limits_{h\rightarrow 0}\frac{f(t+h) - f(t)}{h}
\end{equation}
that uses the present and future values.
We say that the former is causal, while this one is anti-causal.
\end{itemize}
So, when going into a general formulation using time scales as domain, we should
\begin{itemize}
\item use two-sided transforms;
\item define systems with causal derivatives;
\item be careful about the meanings of $s$ and $z$.
\end{itemize}


\subsection{The current Laplace transform on time scales}

The Laplace transform defined on time scales was subject of several publications
but most of them consider the one-sided Laplace transform
\cite{A09,BMT11,BG10,BG10b,BGK11,BP02,D07,DGM10,K11,SR11}.
Anyway, the bilateral was considered in \cite{DGM10b}, where it was defined by
\begin{equation*}
F(s) = \int\limits_{-\infty}^{+\infty} f(t) e_{ \ominus s}^{\sigma}(t) \Delta t,
\end{equation*}
where $e_{s}(t)$ is the generalized exponential
\begin{equation*}
e_{s}(t) = \exp\int\limits_{0}^{t}
\frac{\ln\left( 1+s\mu(\tau)\right) }{\mu(\tau)} \Delta\tau,
\end{equation*}
$e_{ \ominus s}^{\sigma}(t) = e_{ \ominus s}(\sigma(t))$,
and $\ominus s = \frac{-s}{1+s \mu(t)}$.
When the time scale is $\mathbb{R}$, we reobtain \eqref{LTZ_2},
but if the time scale is $\mathbb{Z}$, then we obtain
\begin{equation*}
F(s) = \sum\limits_{n=-\infty}^{+\infty} f(n) (1+s)^{-(n+1)}.
\end{equation*}
Let us make the substitution $z = 1+s$. We have
\begin{equation*}
F(z) = \sum\limits_{n=-\infty}^{+\infty} f(n) z^{-(n+1)}
\end{equation*}
that does not coincide with the $Z$ transform \eqref{LTZ_7}
since we have $n+1$ instead of $n$. This seems to be insignificant
but may create great difficulties as the results we obtain
are different from the conventional ones.


\section{Derivatives and inverses}

In what follows $h\in \mathbb{R}^+$.


\subsection{Nabla and delta derivatives}

We define the nabla derivative by
\begin{equation}
\label{LTZ_15}
f^{\nabla}(t) =
\left\lbrace \begin{array}{cl}
\frac{f(t) - f(\rho(t))}{\nu(t)} & \text{ if } \nu(t) \neq 0\\
\lim\limits_{h\rightarrow 0}\frac{f(t) - f(t-h)}{h}
& \text{ if }  \nu(t) = 0
\end{array}\right.
\end{equation}
and the delta derivative by
\begin{equation}
\label{LTZ_16}
f^{\Delta}(t) =
\left\lbrace \begin{array}{cl}
\frac{f(\sigma(t)) - f(t)}{\mu(t)}
& \text{ if } \mu(t) \neq 0\\
\lim\limits_{h\rightarrow 0}\frac{f(t+h) - f(t)}{h}
& \text{ if } \mu(t) = 0.
\end{array}\right.
\end{equation}
As it can be seen, the first one is causal, while the delta derivative is anti-causal. 
One is the dual of the other \cite{duality}. We can give another form to these derivatives, 
in agreement with our comments done before about time scales. 
We define the nabla derivative by
\begin{equation}
\label{LTZ_17}
f^{\nabla}(t_n)
=
\left\lbrace \begin{array}{cl}
\frac{f(t_n) - f(t_{n-1})}{\nu_n} 
& \text{for any } \;t_n \;\text{that is not left dense},\\
\lim\limits_{h\rightarrow 0^+}\frac{f(t_n) - f(t_n-h)}{h} & \text{for}\; t_n \;
\text{left dense},
\end{array}\right.
\end{equation}
where $\nu_n=t_n-t_{n-1}$, and the delta derivative by
\begin{equation}
\label{LTZ_18}
f^{\Delta}(t_n) =
\left\lbrace \begin{array}{cl}
\frac{f(t_{n+1})-f(t_n)}{\mu_n} &
\text{for any } \;t_n\; \text{that is not right dense},\\
\lim\limits_{h\rightarrow 0^+}\frac{f(t_n+h) - f(t_n)}{h}
& \text{for} \;t_n\; \text{right dense},
\end{array}\right.
\end{equation}
where $\mu_n =t_{n+1}-t_n$. 
The derivatives of higher-order are defined as usual.
Let $f^{\nabla^0} = f^{\Delta^0} = f$. We define
the nabla and delta derivatives of order $n$
as $f^{\nabla^n} = \left(f^{\nabla^{n-1}}\right)^\nabla$
and $f^{\Delta^n} = \left(f^{\Delta^{n-1}}\right)^\Delta$,
$n = 1, 2, \ldots$ In particular, $f^{\nabla^1} = f^{\nabla}$
and $f^{\Delta^1} = f^{\Delta}$.


\subsection{Nabla and delta anti-derivatives}

We are going to obtain the inverses of the above derivatives that we will call
derivatives of order $-1$ or anti-derivatives. To obtain the searched formulae,
we proceed as if $\nu(t)$ and $\mu(t)$ are nonnull. The null situation
will be included later. We rewrite \eqref{LTZ_15} as
\[
\frac{f(t) - f(t-\nu(t))}{\nu(t)} = g(t)
\]
and manipulate it to obtain $f(t) $ as a function of $g(t)$:
\[
f(t) = f(t-\nu(t)) + \nu(t) g(t).
\]
To obtain the nabla and delta inverses, we proceed recursively as stated above.
We then have
\begin{equation*}
f^{\nabla^{-1}}(t)
= \sum\limits_{n=0}^{\infty} \nu_{n+1}(t)f(t - \nu^n(t))
\end{equation*}
and
\begin{equation*}
f^{\Delta^{-1}}(t) = - \sum\limits_{n=0}^{\infty} \mu_{n+1}(t) f(t + \mu^n(t)),
\end{equation*}
where we assume that $f^{\nabla^{-1}}(-\infty) = f^{\Delta^{-1}}(+\infty) = 0$.
It is easy to show that these functions are really the inverses corresponding
to \eqref{LTZ_16} and \eqref{LTZ_17}. To obtain a different formulation,
start from \eqref{LTZ_17} and \eqref{LTZ_18} and proceed similarly. We obtain
\begin{equation*}
f^{\nabla^{-1}}(t_n) = \sum\limits_{m=-\infty}^{n} \nu_mf(t_m)
\end{equation*}
and
\begin{equation}
\label{LTZ_22}
f^{\Delta^{-1}}(t_n) = - \sum\limits_{m=n}^{\infty} \mu_m f(t_m).
\end{equation}
In the constant graininess situation, $\nu^n = \mu^n = nh$,
$\nu_n=\mu_n= h$, and we recover the expressions
obtained in \cite{OrCoTr15} from the generalized fractional derivative.
If $f(t) \equiv 1$, then we conclude immediately that
\[
f^{\nabla^{-1}}(t) = t \qquad \text{and} \qquad
f^{\Delta^{-1}}(t) = -t.
\]
This is a special case that does not verify the above condition: it is nonnull
at infinite. This means that we must be careful when trying to obtain 
the polynomials by computing the anti-derivatives. We return to this subject later.
Those expressions lead us to define the corresponding nabla and delta integrals
over the interval $[a, b]$ respectively by
\begin{equation*}
I^\nabla f = f^{\nabla^{-1}}(b) - f^{\nabla^{-1}}(a)
\end{equation*}
and
\begin{equation*}\
I^\Delta f =f^{\Delta^{-1}}(a) - f^{\Delta^{-1}}(b).
\end{equation*}
Putting $a = b -\nu^k(b)$, one has $\nu_{n+1}(a) = \nu_{n+1+k}(b)$
and $ \nu^{n}(a) = \nu^{n+k}(b)$, and we can write
\[
I^\nabla f =: \int\limits_{\sigma(a)}^{b} f(t)\nabla t 
= \sum\limits_{n=0}^{k-1} \nu_{n+1}(b)f(b - \nu^n(b)).
\]
Similarly, with $ b=a+\mu^k(a) $, we have $ \mu_{n+1}(b) = \mu_{n+1+k}(a)$
and $ \mu^{n}(b) = \nu^{n+k}(a)$ that leads to
\[
I^\Delta f =: \int\limits_{a}^{\rho(b)} f(t)\Delta t 
= \sum\limits_{n=0}^{k-1} \mu_{n+1}(a)f(a + \mu^n(a)).
\]


\section{The nabla and delta general exponentials}
\label{nde}

We consider first the nabla derivative. We want to discover the eigenfunction
of this operator. We begin by introducing the reference instant
$t_0$ and the current instant $t$. The relations between both depend
on the two cases $t > t_0$ or $t < t_0$. We start from $t_0$ and go to $t$
using the values of the graininess. In the first case,
we easily write $t = t_0  + \mu^n(t_0)$, where
$n$ is the number of steps to go from $t_0$ to $t$.
For the second case, we have $t_0-\nu^m(t_0)$.


\subsection{The nabla exponential}
\label{sec:nabla:exp}

Starting from \eqref{LTZ_15}, we write successively
\[
\frac{f(t) - f(t-\nu(t))}{\nu(t)} = s f(t),
\]
\[
f(t) - f(t-\nu(t)) = s\nu(t) f(t),
\]
and
\[
f(t)\left[1 - s\nu(t)\right]  = f(t-\nu(t)).
\]
On the other hand, going back in time, we have
\[
f(t-\nu(t))\left[1 - s\nu(t-\nu(t))\right]  = f(t-\nu(t-\nu(t)))
\]
that can be substituted above to get
\[
f(t)\left[1 - s\nu(t)\right]\left[1 - s\nu(t-\nu(t))\right]
= f(t-\nu(t-\nu(t))).
\]
We can write
\[
f(t) \cdot \prod\limits_{k=1}^{n}\left[1 - s\nu_k(t)\right] 
= f(t-\nu^n(t)).
\]
Consider now a time $t_0$ at which we assume $f(t_0) =1$.
For $t = t_0  + \mu^n(t_0)$ we then have
\[
f(t)=\prod\limits_{k=1}^{n}\left[1 - s\mu_k(t_0)\right]^{-1},
\]
that can be written as
\[
f(t) = \exp \left[ - \sum\limits_{k=1}^{n} \left( 1 - s\mu_k(t_0)\right) \right]
\]
or
\[
f(t) = \exp \left[ - \int\limits_{t_0}^{t}
\frac{\ln\left( 1 - s\mu(\tau)\right)}{\mu(\tau)}\nabla \tau \right].
\]
For $t < t_0$ we must invert the above recursion. In fact, we have
\[
f(t_0)\left[1 - s\nu(t_0)\right]  = f(t_0-\nu(t_0))
\]
and with repetition, and noting that $f(t_0) =1$, we have
\[
\prod\limits_{k=1}^{m}\left[1 - s\nu_k(t_0)\right] = f(t_0-\nu^m(t_0)),
\]
\[
f(t) = \prod\limits_{k=1}^{m}\left[1 - s\nu_k(t_0)\right],
\]
\[
f(t) = \exp \left[ \sum\limits_{k=1}^{m}\ln  \left( 1 - s\nu_k(t_0)\right) \right]
\]
or
\[
f(t) = \exp \left[  \int\limits_{t}^{t_0}
\frac{\ln\left( 1 - s\nu(\tau)\right)}{\nu(\tau)}\nabla \tau \right].
\]
From these equalities, we conclude that the nabla
generalized exponential is given by
\begin{equation}
\label{LTZ_23:a}
e_{\nabla}(t,t_0) =
\begin{cases}
\prod\limits_{k=1}^{n}\left[1 - s\mu_k(t_0)\right]^{-1}
& \text { if } \qquad t= t_n > t_0\\
1 & \text { if } \qquad t = t_0\\
\prod\limits_{k=1}^{m}\left[1 - s\nu_k(t_0)\right]
& \text { if } \qquad t = t_m < t_0
\end{cases}
\end{equation}
or
\begin{equation}
\label{LTZ_23}
e_{\nabla}(t,t_0) =
\begin{cases}
\exp \left[ - \int\limits_{t_0}^{t}
\frac{\ln\left( 1 - s\mu(\tau)\right)}{\mu(\tau)}\nabla \tau \right]
& \text { if } \qquad t = t_n > t_0\\
1 & \text { if } \qquad t = t_0\\
\exp \left[  \int\limits_{t}^{t_0}
\frac{\ln\left( 1 - s\nu(\tau)\right)}{\nu(\tau)}\nabla \tau \right]
& \text { if } \qquad t = t_m < t_0.
\end{cases}
\end{equation}
To illustrate \eqref{LTZ_23:a} and \eqref{LTZ_23},
we are going to compute them for particular cases.

\begin{example}
\label{ex:nabla:exp}
We obtain the nabla generalized exponential 
\eqref{LTZ_23:a}--\eqref{LTZ_23} in different time scales.
\begin{itemize}
\item Let $\mathbb{T} = h \mathbb{Z}$, $h > 0$.
We make $t_0=0$ and put $\nu(t) = \mu(t) = h $, $t=nh$, $n\in \mathbb{Z}$,
leading to $e_{\nabla}(t,0) = \left[1 - sh\right]^{-n}$. Putting $z^{-1}$
for $\left[1 - sh\right]$, as suggested by \eqref{LTZ_7},
we obtain the current discrete-time exponential, $z^n$.

\item Let $\mathbb{T=R}$. Return to the above case and put $h = \frac{t}{n}$. As
$$
\lim\limits_{n\rightarrow \infty}
\left( 1 - \frac{s t}{n} \right)^{-n}  = e^{st},
$$
we obtain $e_{\nabla}(t,0) = e^{st}$.

\item Let $\mathbb{T} = \bigcup_{n=-\infty}^{\infty}
\bigcup_{m=0}^{N} (2n+mh)$ with $N\in \mathbb{N}_0^+$ and $h=1/N$.
Here we also have $t_0=0$ and the graininess functions are given as follows:
\begin{equation}
\label{eq:ex:nu}
\nu(t) =
\begin{cases}
1 & \text{ if } t = 2n, \, n \in \mathbb{Z},\\
h & \text{ otherwise}
\end{cases}
\end{equation}
and
\begin{equation}
\label{eq:ex:mu}
\mu(t) =
\begin{cases}
1 & \text{ if } t = 2n + 1, \, n \in \mathbb{Z},\\
h & \text{ otherwise}.
\end{cases}
\end{equation}
This time scale leads
to the following exponential:
\begin{itemize}
\item if $t= 2n+Mh$, $n \ge 0$, $0 \le M \le N$, then
\[
e_{\nabla}(t,0) = \left[1 - s \right]^{-n} \left[1 - sh \right]^{-nN-M};
\]
\item if $t= -2n+1 - Mh$, $n > 0$, $0 \le M \le N$, then
\[
e_{\nabla}(t,0) = \left[1 - s \right]^{n} \left[1 - sh \right]^{nN+M}.
\]
\end{itemize}

\item Let $\mathbb{T} =  \bigcup_{n=-\infty}^{\infty}  \left[ 2n, 2n+1 \right]$.
This case can be treated from the last one. More precisely, we only have to make
$h$ going  to zero. We start by noting that since $Nh=1$, one has
$t= 2n+Mh = n +nNh+Mh$, which is equivalent to
$nN+M = \frac{t-n}{h}$. We first obtain that
\[
\lim\limits_{h \rightarrow 0} \left[1 - sh \right]^{-\left[\frac{t-n}{h} \right]}
= e^{-s(t-n)}, \quad t\in \left[2n, 2n+1 \right].
\]
Therefore, if $t\ge 0$, then
\[
e_{\nabla}(t,0) = \left[1 - s \right]^{-n} e^{-s(t-n)},
\quad t\in \left[ 2n , 2n+1 \right].
\]
Similarly, if $t< 0$, then
\[
e_{\nabla}(t,0) = \left[1 - s \right]^{n} e^{s(t-n+1)},
\quad t\in \left[ -2n , -2n+1 \right].
\]

\item We now consider a time scale with a periodic graininess $\mu(t+\tau) = \mu(t)$.
Consider a set of positive real numbers $ M = \left\lbrace \mu_i \in \mathbb{R}^+,
i=1, \dots N \right\rbrace$, and the time scale
$$
\mathbb{T}_0 = \left\lbrace \tau_0, \tau_1, \dots \tau_{N-1} \right\rbrace
= \left\lbrace 0, \mu_1, \mu_1+\mu_2, \mu_1+\mu_2+\mu_3,
\ldots, \sum_{i = 1}^{N-1} \mu_i \right\rbrace.
$$
Let $ \tau = \sum_{i=1}^{N} \mu_i$
and prolong the time scale $\mathbb{T}_0$ to left and right
to obtain a time scale $\mathbb{T}$ with periodic graininess $\tau$:
\begin{equation}
\label{eq:ex:period:ts}
\mathbb{T} = \bigcup_{n = -\infty}^{+\infty} \mathbb{T}_0^\tau(n),
\end{equation}
where $\mathbb{T}_0^\tau(n) = \left\{t + n\tau : t \in \mathbb{T}_0 \right\}$.
The nabla exponential in $\mathbb{T}$ is given by:
\begin{itemize}
\item if $t= n\tau + \tau_i$, $0 \le i < N$, $ n=0, 1, \ldots$, then
\[
e_{\nabla}(t,0) = \prod\limits_{k=0}^{N-1}\left[1
- s\mu_k\right]^{-n}\prod\limits_{k=0}^{i}\left[1 - s\mu_k\right]^{-1};
\]
\item if $t= -n\tau + \tau_i$, $0 \le i < N$, $n= 1,2, \ldots$, then
\[
e_{\nabla}(t,0) = \prod\limits_{k=0}^{N-1}\left[1 - s\mu_k\right]^{n}
\prod\limits_{k=0}^{i}\left[1 - s\mu_k\right].
\]
\end{itemize}
\end{itemize}
\end{example}


\subsection{The delta exponential}

Following a similar procedure as in Section~\ref{sec:nabla:exp},
we can obtain the exponential corresponding to the delta derivative.
Starting from \eqref{LTZ_16}, we write successively
\[
f(t+\mu(t)) - f(t)= s \mu(t) f(t)
\]
and
\[
f(t+\mu(t)) = f(t)\left[1 + s\mu(t)\right].
\]
Starting from $t_0$ and repeating the above relation $n$ times, 
we obtain that
\[
f(t_0+\mu^n(t_0)) =  f(t_0)
\cdot \prod\limits_{k=1}^{n}\left[1 + s\mu_k(t_0)\right],
\]
leading to
\[
f(t) = \prod\limits_{k=1}^{n}\left[1 + s\mu_k(t_0)\right].
\]
Similarly to the nabla case, it is not a hard task to conclude that
if $t < t_0$, then
\[
f(t) = \prod\limits_{k=1}^{n}\left[1 + s\nu_k(t_0)\right]^{-1},
\]
leading to the delta exponential
\begin{equation}
\label{LTZ_24}
e_{\Delta}(t,t_0)
=  \begin{cases}
\prod\limits_{k=1}^{n}\left[1 + s\mu_k(t_0)\right]
& \text{ if } \qquad t= t_n > t_0\\
1 & \text{ if } \qquad t = t_0\\
\prod\limits_{k=1}^{m}\left[1 + s\nu_k(t_0)\right]^{-1}
& \text{ if } \qquad t= t_m < t_0
\end{cases}
\end{equation}
or
\begin{equation}
\label{LTZ_25}
e_{\Delta}(t,t_0)
=  \begin{cases}
\exp \left[ \int\limits_{t_0}^{t}
\frac{\ln\left( 1 + s\mu(\tau)\right)}{\mu(\tau)}\Delta \tau \right]
& \text{ if } \qquad t= t_n > t_0\\
1 & \text{ if } \qquad t = t_0\\
\exp \left[ - \int\limits_{t}^{t_0}
\frac{\ln\left( 1 + s\nu(\tau)\right)}{\nu(\tau)}\Delta \tau \right]
& \text{ if } \qquad t= t_m < t_0.
\end{cases}
\end{equation}

Let us see what happens in the particular cases
considered in Example~\ref{ex:nabla:exp}.

\begin{example}
\label{ex:delta:exp}
We obtain the delta exponential \eqref{LTZ_24}--\eqref{LTZ_25}
in different time scales.
\begin{itemize}
\item If $\mathbb{T}=h\mathbb{Z}$, $h > 0$, then we make again
$t_0=0$, $ \nu(t) = \mu(t) = h $, and $t=nh$, $n\in \mathbb{Z}$,
leading to $e_{\Delta}(t,0) = \left[1 + sh\right]^{n}$.
Putting $z$ for $ \left[1 + sh\right]$, as suggested by \eqref{LTZ_10},
we obtain once again the current discrete-time exponential, $z^n$.

\item If $\mathbb{T}=\mathbb{R}$, then
we return to the previous case and put $h = \frac{t}{n}$. As
$$
\lim\limits_{n\rightarrow \infty}
\left[ 1 + \frac{s t}{n} \right]^{n}  = e^{st},
$$
we obtain that $e_{\Delta}(t,0) = e^{st}$.

\item Let $\mathbb{T} = \bigcup_{n=-\infty}^{\infty} \bigcup_{m=0}^{N} 2n+mh$
with $N\in \mathbb{N}_0^+$ and $h=1/N$, and $t_0=0$.
The graininess functions $\nu$ and $\mu$ are given by
\eqref{eq:ex:nu} and \eqref{eq:ex:mu}, respectively.
This time scale leads to the following delta exponential:
\begin{itemize}
\item if $t= 2n+Mh$, $n \ge 0$, $0 \le M \le N$, then
\[
e_{\Delta}(t,0) = \left[1 + s \right]^{n} \left[1 + sh \right]^{nN+M};
\]
\item if $t= -2n+1 - Mh$, $n > 0$, $0 \le M \le N$, then
\[
e_{\Delta}(t,0) = \left[1 + s \right]^{n} \left[1 + sh \right]^{-nN-M}.
\]
\end{itemize}

\item Let $\mathbb{T} =  \bigcup_{n=-\infty}^{\infty}  \left[ 2n, 2n+1 \right]$.
This case can be treated from the last one. We only have to make $h$ tend to zero.
We start by noting that $ Nh=1 $ and $t= 2n+Mh = n +nNh+Mh$, that is,
$nN+M = \frac{t-n}{h}$. We first obtain that
\[
\lim\limits_{h \rightarrow 0} \left[1 + sh \right]^{\left[\frac{t-n}{h} \right]}
= e^{s(t-n)}, \quad t\in \left[2n, 2n+1 \right],
\]
concluding that if $t\ge 0$, then
\[
e_{\Delta}(t,0) = \left[1 + s \right]^{n} e^{s(t-n)},
\quad t\in \left[ 2n , 2n+1 \right].
\]
Similarly, if $t< 0$, then
\[
e_{\Delta}(t,0) = \left[1 + s \right]^{-n} e^{-s(t-n+1)},
\quad t\in \left[ -2n , -2n+1 \right].
\]
\item Consider now the time scale \eqref{eq:ex:period:ts}
with periodic graininess $\mu(t+\tau) = \mu(t)$.
The corresponding delta exponential is given by:
\begin{itemize}
\item if $t= n\tau + \tau_i$, $0 \le i < N$, $ n=0, 1, \ldots$, then
\[
e_{\Delta}(t,0) = \prod\limits_{k=0}^{N-1}\left[1 + s\mu_k\right]^{n}
\prod\limits_{k=0}^{i}\left[1 + s\mu_k\right];
\]
\item if $t= -n\tau + \tau_i$, $0 \le i < N$, $n= 1,2, \ldots$, then
\[
e_{\Delta}(t,0) = \prod\limits_{k=0}^{N-1}\left[1 + s\mu_k\right]^{-n}
\prod\limits_{k=0}^{i}\left[1 + s\mu_k\right]^{-1}.
\]
\end{itemize}
\end{itemize}
\end{example}


\subsection{Properties of the exponentials}
\label{pe}

The exponentials \eqref{LTZ_23:a}--\eqref{LTZ_23} and
\eqref{LTZ_24}--\eqref{LTZ_25} have interesting properties
that we describe in the following. To enhance
the importance of the parameter $s$, in the sequel we use 
the notations $e_\nabla (t,t_0;s)$ and $e_\Delta (t,t_0;s)$.
\begin{enumerate}
\item \textit{Changing the role of instants}.
If we change the role of $t$ and $t_0$, then
\begin{equation*}
e_{\nabla}(t_0,t;s) = 1/e_{\nabla}(t,t_0;s).
\end{equation*}
The delta exponential satisfies similar relation:
$e_{\Delta}(t_0,t;s) = 1/e_{\Delta}(t,t_0;s)$.

\item \textit{Relation between nabla and delta exponentials}.
The function defined in \eqref{LTZ_24} with the substitution of
$-s$ for $s$ is the inverse of \eqref{LTZ_22}:
\begin{equation*}
e_{\Delta}(t,t_0;s) = 1/e_{\nabla}(t,t_0;-s).
\end{equation*}
\item As $h \rightarrow 0$,
both nabla and delta exponentials converge to $e^{st}$.
This was already seen in Examples~\ref{ex:nabla:exp}
and \ref{ex:delta:exp}.
\item Attending to the way how both exponentials
were obtained, we can easily conclude that
\begin{equation*}
\left[e_{\nabla}(t,t_0;s)\right]^{\nabla^{N}}
= s^N e_{\nabla}(t,t_0;s)
\end{equation*}
and
\begin{equation*}
\left[e_{\Delta}(t,t_0;s)\right]^{\Delta^{N}}
= s^N e_{\Delta}(t,t_0;s),
\end{equation*}
where $\nabla^{N}$ and $\Delta^{N}$ denote respectively
the $N$th order nabla and delta derivative,
$N \in \mathbb{N}$. Later we try to generalize this property.
We can say that the exponentials are also eigenfunctions
of the higher-order derivatives.

\item \textit{Translation}.
An exponential defined on a given time scale
is called translation of $e_{\nabla}(t,t_0;s)$
if it is obtained by a similar procedure but using a different 
reference instant $\tau_0\in \mathbb{T}$. This translation 
can be obtained by a product of exponentials (cf. \eqref{LTZ_29}).

\item It is important to know how exponentials increase/decrease
for $s \in \mathbb{C} $. Let us consider the nabla case.
The delta situation is similar and follows straightforwardly,
so we omit it here.  Let $h_M = \max(\nu_k, \mu_l)$ 
and $h_m = \min(\nu_k, \mu_l)$, $k,l \in \mathbb{Z}$. 
The nabla exponential $e_{\nabla}(t,t_0;s)$
\begin{itemize}
\item is a real number for any real $s$;
\item is positive for any real number $s$ such that $s < \frac{1}{h_M}$;
\item oscillates for any real number $s$ such that $s > \frac{1}{h_m}$;
\item is bounded for values of $s$ inside the inner Hilger circle
$\left| 1-sh_m\right|  =1$;
\item has an absolute value that increases as $|s|$ increases
outside the outer Hilger circle $\left| 1-sh_M\right| =1$,
going to infinite as $|s| \rightarrow \infty$.
\end{itemize}
This means that each graininess defines a circle centred at its inverse
and passing at $ s=0$. In general, we have infinite circles that reduce
to one when the graininess is constant. When it is null, it degenerates
into the imaginary axis.

\item \textit{Product of exponentials}.
The following relations hold:
\begin{equation}
\label{LTZ_29}
\begin{split}
e_{\nabla}(t,t_0;s) \cdot e_{\Delta}(\tau,t_0;-s)
&= e_{\nabla}(t,\tau;s)\\
&= e_{\Delta}(\tau,t;-s).
\end{split}
\end{equation}
Note that the product $e_{\nabla}(t,t_0;s) \cdot e_{\Delta}(t,t_0;v)$, $s\ne v$,
is well defined but does not lead, in general, to an exponential as we defined
above, because its parameter $s$ (eigenvalue) is time dependent. The same
happens with the product of exponentials of the same type:
$e_{\nabla}(t,t_0;s) \cdot e_{\nabla}(t,t_0;v)$ or
$e_{\Delta}(t,t_0;s) \cdot e_{\Delta}(t,t_0;v)$.
Now let $t_n > t_m$. We have
\begin{equation}
\label{LTZ_30}
\begin{split}
e_{\nabla}(t_n,t_0;s)
&= e_{\nabla}(t_m,t_0;s) \cdot e_{\nabla}(t_n,t_m;s)\\
&= e_{\nabla}(t_m,t_0;s) \cdot e_{\Delta}(t_m,t_n;-s)
\end{split}
\end{equation}
and
\begin{equation}
\label{LTZ_31}
\begin{split}
e_{\Delta}(t_n,t_0;s)
&= e_{\Delta}(t_m,t_0;s) \cdot e_{\Delta}(t_n,t_m;s)\\
&= e_{\Delta}(t_m,t_0;s) \cdot e_{\nabla}(t_m,t_n;-s).
\end{split}
\end{equation}
In particular, if $t_m = t_n - \tau$, then
\[
e_{\nabla}(t_n,t_0;s)
= e_{\nabla}(t_n - \tau,t_0;s) \cdot e_{\nabla}(t_n,t_n - \tau;s)
\]
and
\begin{equation}
\label{LTZ_32}
e_{\nabla}(t_n - \tau,t_0;s)
= e_{\nabla}(t_n,t_0;s) \cdot e_{\nabla}(t_n,t_n - \tau;-s).
\end{equation}
Relation \eqref{LTZ_32} states the \emph{translation or shift property}
of an exponential, and will be used later.

\item \textit{Scale change}.
Let $ a $ be a positive real number. We have:
\[
e_{\nabla}(a t,t_0;s) = e_{\nabla}(t,t_0;a s).
\]
If $a < 0$, a similar relation is obtained
but involving the delta exponential.

\item \textit{Cisoids}. For $\mathbb{T=R}$
the exponentials degenerate into the cisoids $e^{i\omega t}$
when $s=i\omega$. Similarly, in $\mathbb{T}= h \mathbb{Z}$,
$h > 0$, we obtain  $e^{i\omega nh}$.
We put the question: what happens in the general case?
The answer is easy: \emph{on a general time scale there are no cisoids}.
Thus, in general, there are no periodic signals. This has a very important 
consequence: \emph{on a general time scale we cannot define Fourier transforms},
which is a powerful drawback of the theory.
\end{enumerate}

The properties of the exponentials are here enlarged for measure
chains $\mathbb{T}$. In the real line case, i.e., 
$\mathbb{T} = \mathbb{R}$, we obtain the ordinary exponential, whose
properties are well known. Very interesting and important is also the
case $\mathbb{T} = h \mathbb{Z}$, studied in \cite{OrCoTr15}.


\section{Suitable general Laplace transforms}
\label{sec:gLt}

To define a Laplace transform, we adopt a point of view somehow different
from the usual. We obtained exponentials that are eigenfunctions of the nabla
and delta derivatives. This means that we would like to express a given function
in terms of these exponentials. To do it, we start from the inverse transform.


\subsection{The nabla inverse transform}
\label{nit}

Assume that we have a signal, $f(t)$, defined on a given
time scale, that has a given transform $F_\nabla(s)$ (undefined for now;
defined later by \eqref{eq:nablaLT}). Its inverse transform would serve
to synthesize it from a continuous set of elemental exponentials:
\begin{equation}
\label{LTZ_33}
f(t) = \frac{1}{2\pi i} \oint_C F_\nabla(s) e_\nabla(t,t_0;s) ds,
\end{equation}
where the integration path, $C$, is a simple closed contour in a region
of analyticity of the integrand and the poles introduced by the exponential.
Consider the Hilger circle with radius 
$R > \frac{1}{h_m} = \frac{1}{\min(\nu_k, \mu_l)}$,
centred at $\frac{1}{R}$, and described in the counterclockwise direction.
We will not consider other circles here, leaving the subject for future studies.
We call this region the \emph{fundamental region}, and we represent it by $R_c$.
In the continuous case it degenerates into a vertical straight line. Let us start
with the simple situation corresponding to $F_\nabla(s) \equiv 1$.
We are expecting to obtain an impulse as the inverse transformed function.
Let us try to compute it. As it is easy to observe,
the integral is zero for $t \leq t_0$, because the integrand
is analytical. For $t > t_0$, the integrand has $n$ poles at the points
$s_k = \frac{1}{\nu_k(t_0)}$, $k=1, 2, \ldots$ If $n > 1$ and the
integration path includes all the poles, then the integral is zero. So, only for
$n = 1$ the integral is nonnull. This means that to have a pole at $ n=0$ we have
to introduce one. This leads to change the integrand in the above formula
\eqref{LTZ_33} into
\begin{equation}
\label{LTZ_34}
f(t) = -\frac{1}{2\pi i } \oint_C F_\nabla(s) e_\nabla(t+\mu(t),t_0;s) ds,
\end{equation}
which is our definition of inverse nabla Laplace transform.
It is a simple task to show that the inverse nabla Laplace transform
of $F_\nabla(s) \equiv 1$ is the impulse $\delta(t_n-t_0) = \frac{1}{\nu_1(t_0)}$.
We can show that it is the derivative of the Heaviside unit step defined by
\begin{equation*}
\epsilon(t,t_0)
= \begin{cases}
1 & \text{ if } \;\;\; t \ge t_0,\\
0 & \text{ if } \;\;\; t < t_0,
\end{cases}
\end{equation*}
$t \in \mathbb{T}$. Now we are going to compute the nabla derivative
of $e_\nabla(t+\mu(t))$. We can immediately write that
\[
\left[ e_\nabla(t+\mu(t))\right]^\nabla
= \frac{e_\nabla(t+\mu(t),t_0)
- e_\nabla(t+\mu(t)-\nu(t+\mu(t)),t_0)}{\nu(t+\mu(t))}.
\]
Let $ t \geq t_0 $. We have
\[
\left[e_\nabla(t+\mu(t))\right]^\nabla
= \frac{\prod\limits_{k=1}^{n+1}\left[1 - s\mu_k(t_0)\right]^{-1}
- \prod\limits_{k=1}^{n}\left[1 - s\mu_k(t_0)\right]^{-1}}{\nu_{n+1}(t_0)}
= s\cdot e_\nabla(t+\mu(t)).
\]
For $ t < t_0 $ the situation is similar. With this result we prove one of the
most important properties of the nabla Laplace transform $NLT$:
\begin{equation}
\label{LTZ_35}
NLT\left[f^\nabla(t) \right] = s F_\nabla(s),
\end{equation}
where $F_\nabla(s)$ is the nabla Laplace transform of $f(t)$, i.e.,
$F_\nabla(s) = NLT[f(t)]$.


\subsection{Properties of the nabla Laplace transform}

We deduce the properties of the nabla Laplace transform
from the inversion integral. Precisely, the following properties hold 
(we assume $s$ to be inside the region of convergence of the nabla Laplace
transform; this region will be defined soon):
\begin{itemize}
\item \textit{Linearity}.
It is an evident property that does not need any computation.
\item \textit{Transform of the derivative}.
Attending to the equality \eqref{LTZ_35},
we easily deduce, by repeated application of \eqref{LTZ_35},
that if $N\in \mathbb{N}_0$, then
\begin{equation*}
NLT\left[f^{\nabla^{N}}(t)\right] = s^N F_\nabla(s),
\end{equation*}
restating a well known result in the context of the Laplace transform.

\item \textit{Time scaling}.
Let $ a $ be a positive real number.\footnote{It could be negative,
but it is not interesting.} We have
\begin{equation*}
NLT\left[f(a t)\right] = \frac{1}{a} F_\nabla\left(\frac{s}{a}\right).
\end{equation*}
\end{itemize}


\subsection{The nabla direct transform}

The shift property \eqref{LTZ_32} gives a suggestion of how to define
the nabla Laplace transform. In fact it must be defined in terms
of the delta exponential, due to property \eqref{LTZ_29}
particularised to $t \rightarrow t+\mu(t) \: \text{and} \: \tau = t$.  
We then have
\begin{equation}
\label{eq:nablaLT}
F_\nabla(s) = \sum\limits_{n=-\infty}^{+\infty}
\nu_n f(t_n) e_\Delta(t_n,t_0;-s).
\end{equation}
The region in the complex plane where the series \eqref{eq:nablaLT} converges,
is called ``region of convergence'' of the nabla Laplace
transform. In a general time scale, as we are assuming,
we cannot obtain current properties of the Laplace or $Z$ transforms.
Particularising for the time scale $h \mathbb{Z}$ we easily deduce
the \textit{modulation or shift in} $s$ property \cite{OrCoTr15}:
\begin{equation*}
NLT\left[f(nh)e_{\nabla}(nh,-s_0) \right]
= F_\nabla\left(\frac{s-s_0}{1-s_0h}\right).
\end{equation*}
To obtain this result we only have to use the inverse 
transform \eqref{LTZ_30}. Similarly, we can obtain that
\[
NLT\left[f(nh)e_{\Delta}(nh,-s_0) \right]
= F_\nabla(s+s_0-ss_0h).
\]
Concerning the \textit{product or convolution in} $s$, we easily obtain that
\begin{equation*}
NLT\left[f(nh)g(nh) \right] = - \frac{1}{2\pi i}
\oint F_\nabla(v) G_\nabla\left(\frac{v-s}{1-vh}\right)
\frac{1}{1-vh} dv.
\end{equation*}
To obtain this result, its enough to insert \eqref{LTZ_34}
into the first member and use the second equality of \eqref{LTZ_30}.


\subsection{Examples}
\label{exn}

In general, it is difficult to compute the direct transform.
So we are going to invert the problem: we look
for the inverses of well known functions.

\begin{itemize}
\item \textit{Causal and anti-causal exponentials}.
A nabla causal exponential is defined by the inverse transform of the function
$\frac{1}{s-p}$. Begin by assuming that $p$ is in the region outside
the Hilger circle defined in Section~\ref{nit}. We have:
\[
f(t) = -\frac{1}{2\pi i } \oint_C \frac{1}{s-p} e_\nabla(t+\mu(t),t_0;s) ds.
\]
If $t < t_0$, then there are no poles in the region delimited by the integration
path $C$ and the function is identically null.
If $t \ge t_0$, then we can have several poles inside and one outside the region
bounded by the integration path. The residue corresponding to $p$ is equal
to the negative sum of the residues for the inner poles.
We conclude that the nabla Laplace transform $NLT$ of the causal exponential
is given by
\begin{equation}
\label{LTZ_41}
NLT \left[ e_\nabla(t,t_0;p)\epsilon(t,t_0) \right] = \frac{1}{s-p}
\end{equation}
for $ s $ inside $R_c$. In fact, and as a simple reasoning can show,
the region of convergence is larger than $R_c$, since the integrand is analytical
outside the $C$ circle. This can be enlarged until it meets the point $s=p$.
We can say that the region of convergence $R_p$ is the interior of a circle
with center at $r_p = \frac{|p|^2}{2\Re(p)}$ and passes on $p$.
With a likewise reasoning, we obtain a similar result for the anti-causal
exponential. In this case the pole must be inside the region referred above.
We conclude that when $t \ge t_0$ all the poles of the integrand are inside
the integration region. The integral is zero for $t< t_0$ and there is only
one pole $s=p$ in this case. We obtain
\begin{equation*}
NLT \left[- e_\nabla(t,t_0;p)\epsilon(-t-\mu(t),t_0) \right] = \frac{1}{s-p}
\end{equation*}
for $s$ outside $R_p$.

\item \textit{Unit steps}.
From the above results, a passage to the limit as $p \rightarrow 0$ allow us
to conclude that the nabla Laplace transform of the unit step is $1/s$:
\begin{equation*}
NLT \left[ \epsilon(t) \right] = \frac{1}{s}
\end{equation*}
for $ s $ inside $R_c$ and
\begin{equation*}
NLT \left[- \epsilon(-t-\mu(t)) \right] = \frac{1}{s}
\end{equation*}
for $ s $ outside $R_c$.
\item \textit{Causal power-exponential}.
Take \eqref{LTZ_41} and compute the usual derivative with respect to $p$.
The right hand side is equal to $\frac{1}{(s-p)^2}$. The derivative of the
exponential in the left hand side can be computed with the help of the logarithm:
\[
\frac{d}{d p} \ln e_\nabla(t,t_0;p) = \sum_{k=0}^n \mu_k,
\]
\[
\frac{d}{d p}e_\nabla(t,t_0;p)
= \sum_{k=0}^n \mu_k e_\nabla(t,t_0;p)
= (t-t_0) e_\nabla(t,t_0;p)
\]
and
\begin{equation*}
NLT \left[ (t-t_0)e_\nabla(t,t_0;p)\epsilon(t,t_0) \right] = \frac{1}{(s-p)^2}
\end{equation*}
for $ s $ inside $R_c$. We can generalise this result for any order
and also for anti-causal exponentials.
\end{itemize}
Extending the approach of \cite{Orti:14}, these results can be used to compute
the impulse response of a linear system described by
differential equations of constant coefficients.
Consider a dynamic equation on time scales
with the general format
\begin{equation*}
\sum\limits_{k=0}^{N} a_k \left(y(t)\right)^{\nabla^k}
= \sum\limits_{k=0}^{M} b_k \left(x(t)\right)^{\nabla^k},
\end{equation*}
where $a_N = 1$. The orders $N$ and $M$ are any given positive integers.
With the nabla Laplace transform, we obtain the transfer function
\begin{equation}
\label{LTZ_47}
H(s) = \frac{Y(s)}{X(s)}
= \frac{\sum\limits_{k=0}^{M} b_k s^{k}}{\sum\limits_{k=0}^{N} a_k s^{k}}.
\end{equation}
Its inverse gives the impulse response of the system \cite{OrCoTr15,Rob03}.
To obtain it, we fix a region of convergence, decompose it in partial fractions,
and invert each one using the above results. Most of the results presented
in \cite{OrCoTr15} apply also here and are left to the reader.


\subsection{Uniqueness of the nabla transform}

From relation \eqref{LTZ_34}, we state a very important conclusion:
\begin{quotation}
\textit{For a given nabla Laplace transform, there are doubly infinite inverse
transforms, accordingly to the chosen time scale and the region of convergence.}
\end{quotation}
This means that with a given region of convergence and a fixed time scale,
the inverse nabla Laplace transform is unique. Changing the region of convergence,
we obtain a function with different characteristics.
Let us return back to the shift property \eqref{LTZ_32} and rewrite it in the format
\begin{equation*}
NLT\left[f(\bar{t} ) \right]
= e_{\Delta}(\bar{t}+\tau,\bar{t};-s) F_\nabla(s).
\end{equation*}
This shows how to interpolate the original function $f(t)$, $t\in \mathbb{T}$,
to a new time scale defined by the \textquotedblleft differences\textquotedblright
\ between all the elements in $\mathbb{T}$. Let $\overline{\mathbb{T}}$ be the
super time scale associated with $\mathbb{T}$, i.e.,
$\overline{\mathbb{T}} = \left\lbrace \bar{t} :
\bar{t} = t - \tau \text{ with } t, \tau \in\mathbb{T} \right\rbrace $.
This means that we have a function $ \bar{f}(t) $ defined in
$\overline{\mathbb{T}}$ such that $\bar{f}(t) = f(t)  \: \text{for} \: t\in\mathbb{T}$,
because $\mathbb{T} \subset \overline{\mathbb{T}}$. This function can be considered
as an interpolated version of the original function.
Each new time grid  extracted from $\overline{\mathbb{T}}$ and equivalent to
$\mathbb{T}$ originates another function with the same transform.
We can write:
\[
\sum\limits_{n=-\infty}^{+\infty} \nu_n f(t_n) e_\Delta(t_n,t_0;-s)
= \sum\limits_{n=-\infty}^{+\infty} \bar{\nu}_n f(\tau_n) e_\Delta(\tau_n,t_0;-s).
\]
Multiplying both sides by $e_\nabla(t_{n+1},t_0;s)$ and integrating accordingly
to \eqref{LTZ_34}, we obtain that
\[
f(t_n) = -\sum\limits_{k=-\infty}^{+\infty} \bar{\nu}_k f(\tau_k)
\frac{1}{2\pi i } \oint_C  e_\Delta(\tau_k,t_0;-s) e_\nabla(t_{n+1},t_0;s) ds
\]
and
\begin{equation*}
f(t_n) = -\sum\limits_{k=-\infty}^{+\infty} \bar{\nu}_k f(\tau_k)
\frac{1}{2\pi i } \oint_C  e_\nabla(t_{n+1},\tau_k;s) ds.
\end{equation*}
In particular, for $t_n = nh$ (we can choose $h$ as the average
of $\bar{\nu}_k$)
\begin{equation*}
f(nh) = -\sum\limits_{k=-\infty}^{+\infty} \bar{\nu}_k f(\tau_k)
\frac{1}{2\pi i } \oint_C  e_\nabla((n+1)h,\tau_k;s) ds.
\end{equation*}
This shows how we can convert a nonuniformly sampled function into another
regularly sampled (cf. \cite{Orti01}). This is very important, because in this
framework, we can use the results presented in \cite{OrCoTr15}. The function
\begin{equation*}
\varphi((n+1)h,\tau_k)  = \bar{\nu}_k  \frac{1}{2\pi i }
\oint_C  e_\nabla((n+1)h,\tau_k;s) ds
\end{equation*}
is a generalisation of the well-known sinc function and \eqref{LTZ_47}
can be considered as a generalisation
of the Shanon--Whitacker theorem \cite{Rob03}.


\subsection{Existence}

As for both the standard Laplace and $Z$ transforms,
also for the nabla Laplace transform it is not an easy task
to formulate necessary and sufficient conditions of existence.
We can easily state, however, sufficient conditions.
Let us consider a function $f(t)$ assuming finite values
for every finite $t_n$, $n\in \mathbb{Z}$. We assume that
it is asymptotically bounded by two exponentials:
\[
\left| f(t) \right| < A \cdot e_\nabla(t,t_0;a), \quad t>t_{n_a},
\]
\[
\left| f(t) \right| < B \cdot e_\nabla(t,t_0;b), \quad t<t_{n_b}.
\]
We can write that
\begin{equation*}
\begin{split}
\Biggl|  \sum\limits_{n=-\infty}^{+\infty} & f(t_n) e_\Delta(t_n,t_0;-s) \Biggr|
\leq  \sum\limits_{n=-\infty}^{+\infty} \left|  f(t_n) e_\Delta(t_n,t_0;-s) \right|\\
&< \sum\limits_{n=n_b}^{n_a} \left|  f(t_n) e_\Delta(t_n,t_0;-s) \right|
+ A \sum\limits_{n=-\infty}^{n_b-1} \left|  e_\nabla(t,t_0;b) e_\Delta(t_n,t_0;-s) \right|\\
&\qquad + B \sum\limits_{n=n_a+1}^{+\infty} \left|e_\nabla(t,t_0;a) 
e_\Delta(t_n,t_0;-s) \right|.
\end{split}
\end{equation*}
According to the results obtained in Sections~\ref{pe} and \ref{sec:gLt},
we conclude that the two infinite summations converge inside and outside two
circular regions with center on the positive real line and passing at $s=0$.

\begin{remark}
As pointed out in \cite{Bas12,D07}, there are two definitions
of ``exponential function'': ``type I'' and ``type II''. The first
is essentially the definition currently used when the time scale is 
$\mathbb{R}$, which does not make sense on a general time scale. 
Only the type II, as we do here, is suitable.
\end{remark}


\subsection{The delta transform}

We studied carefully the transform suitable for dealing with causal systems.
This does not mean that we cannot formulate another transform taking the
delta exponential as the basis. It is not very difficult to realise its formulation.
Similarly as before, each function $f(t)$ can be synthesized by
\begin{equation*}
f(t) =  \frac{1}{2\pi i } \oint_C F_\Delta(s) e_\Delta(t-\nu(t),t_0;s) ds,
\end{equation*}
where the integration path is the one used for the nabla transform.
The relation between the two exponentials
suggests that we define the direct transform by
\begin{equation}
\label{LTZ_52}
F_\Delta(s) = \sum\limits_{n=-\infty}^{+\infty}
\nu_n f(t_n) e_\nabla(t_n,t_0;-s).
\end{equation}


\subsection{Shift, convolution, and correlation}

Although we already discussed something about the shift property,
we need to tell something more. Consider $f(t)$
as an application from $\mathbb{R}$ to $\mathbb{R}$.
When we write  $f(t-t_0)$, $t_0\in \mathbb{R}$, we know exactly its meaning
and we can associate a geometrical interpretation to it. We have no problem
because both $t-t_0$ and $f(t-t_0)$ are well defined. Now consider the situation
we are dealing with: $f(t)$ is an application from $\mathbb{T}$ to $\mathbb{R}$.
In this case it may happen that $t-t_0$ does not belong to $\mathbb{T}$
and consequently $f(t-t_0)$ is not defined. This means that we must devise a meaning
for $f(t-t_0)$. As we mentioned before, mainly in Section~\ref{nde},
we took a reference point $t_0$ that we used to state any instant
in terms of the graininess. So, instead of $f(t)$, we should write $f(t-t_0)$.
When writing  $f(t-t_5)$, for example, we are meaning that a new reference point
was took. We can say that we produced a shift. In conclusion, a shift is nothing
else than a change in the reference point. As consequence we can say that
\[
F_\nabla(s) = \sum\limits_{n=-\infty}^{+\infty} \nu_n f(t_n) e_\Delta(t_n,t_0;-s)
=  \sum\limits_{n=-\infty}^{+\infty} \nu_n f(t_n-t_N) e_\Delta(t_n,t_N;-s).
\]
Our main problem here is to compute the nabla Laplace transform of a shifted function.
Let us consider a function $g(t_n) = f(t_n-t_m)$ for some instants 
$t_n, t_m \in \mathbb{T}$. The nabla Laplace transform of $g(t)$ is
\[
G_\nabla(s) = \sum\limits_{n=-\infty}^{+\infty} \nu_n f(t_n-t_m) e_\Delta(t_n,t_0;-s).
\]
From \eqref{LTZ_31} we can write
\[
e_\Delta(t_n,t_0;-s)  = e_\Delta(t_m,t_0;-s) e_\Delta(t_n,t_m;-s)
\]
that leads to
\begin{equation}
\label{LTZ_53}
\begin{split}
G_\nabla(s) &= \sum\limits_{n=-\infty}^{+\infty} \nu_n f(t_n-t_m)
e_\Delta(t_n,t_m;-s) e_\Delta(t_m,t_0;s) \\
&=  e_\Delta(t_m,t_0;s) F_\nabla(s).
\end{split}
\end{equation}
So, as in the classical case, the shift in time produces a multiplication
by an exponential in the transform domain.

Attending to \eqref{LTZ_32} and \eqref{LTZ_34}, we have
\[
f(t - t_m )  =  -\frac{1}{2\pi i } \oint_C F_\nabla(s)
e_{\nabla}(t+\mu(t),t_m;s) e_\Delta(t_m,t_0;s)  ds.
\]
Inserting the nabla Laplace transform expression inside the integral,
we obtain that
\begin{equation*}
f(t - t_m ) = \sum\limits_{n=-\infty}^{+\infty}\nu_n f(t_n)
-\frac{1}{2\pi i }\oint_C e_\Delta(t_n,t_0;-s)
e_{\nabla}(t+\mu(t),t_m;s) e_\Delta(t_m,t_0;s)ds,
\end{equation*}
where we attend to the fact that the nabla Laplace transform converges
uniformly in the region of convergence. Then,
\begin{equation}
\label{LTZ_54}
f(t - t_m )  = \sum\limits_{n=-\infty}^{+\infty} \nu_n f(t_n) \phi(t,t_m)
\end{equation}
with
\[
\phi(t,t_m) = -\frac{1}{2\pi i } \oint_C
e_\Delta(t_n,t_0;-s) e_{\nabla}(t+\mu(t),t_m;s) e_\Delta(t_m,t_0;s) ds.
\]
Using \eqref{LTZ_31},
\begin{equation*}
\begin{split}
e_\Delta(t_n,t_0;-s) e_\Delta(t_m,t_0;s)
&= e_\Delta(t_n,t_0;-s) e_\Delta(t_n,t_0;s) e_\Delta(t_m,t_n;-s)\\
&= e_\Delta(t_n,t_0;-s^2) e_\Delta(t_m,t_n;-s),
\end{split}
\end{equation*}
\begin{equation*}
\phi(t,t_m) = -\frac{1}{2\pi i } \oint_C   e_{\nabla}(\sigma (t_n),t_m;s)
e_\Delta(t_m,t_n;-s) e_\Delta(t_n,t_0;-s^2) ds.
\end{equation*}
As we can see, \eqref{LTZ_54} is an interpolation formula that obtains
function values for the extended time scale in a way that keeps constant
the nabla Laplace transform. The function $\phi(t,t_m)$
has the role of an interpolating function.

Now it is a simple task to define convolution $f(t) \ast g(t)$ and correlation
$f(t) \ast g(-t)$ on time scales, having in mind to preserve the property
\begin{equation*}
NLT\left[f(t) \ast g(t) \right] = F_\nabla(s) G_\nabla(s).
\end{equation*}
We define convolution by
\begin{equation*}
f(t) \ast g(t)  = \sum\limits_{m=-\infty}^{+\infty}  \nu_m  f(t_m) g(t_n - t_m),
\end{equation*}
where $g(t_n - t_m)$ is defined by \eqref{LTZ_54}. With this convolution
we can show that the output of the system defined by \eqref{LTZ_52}
is given, as usually, by the convolution of the input and the
impulse response defined as the inverse of \eqref{LTZ_53}.
For the correlation the situation is similar but we have to involve
both transforms \cite{OrCoTr15}. In fact we can define it by the convolution
with the reflected function $g(-t)$:
\begin{equation*}
NLT\left[f(t_k) \ast g(-t_k) \right]
= F_\nabla(s) G_\Delta(-s).
\end{equation*}


\subsection{Higher-order and fractional derivatives}

The study of fractional derivatives on arbitrary time scales
is a subject of great interest \cite{BMT11,MyID:296,Bal14,WuBaZe}.
Having defined convolution, we are able to use it to obtain a general formula
for nabla and delta higher-order and fractional derivatives through
the inversion of $s^N$, $N\in \mathbb{Z}_0^+$, and in general
$s^\alpha$, $\alpha\in \mathbb{R}$. As the first case is a particular case,
we solve the second that is more general. Let $\delta_\nabla^{(\alpha)}(t)$
be the causal inverse of such function in order to obtain the nabla derivative.
With the anti-causal we would obtain the delta derivative. We have
\begin{equation*}
\delta_\nabla^{(\alpha)} (t)
= -\frac{1}{2\pi i } \oint_\gamma s^\alpha e_\nabla(t+\mu(t),t_0;s) ds,
\end{equation*}
allowing us to define a general fractional derivative of order $\alpha$ by
\begin{equation*}
x^{(\alpha)}(t) = \delta_\nabla^{(\alpha)} (t) \ast x(t),
\end{equation*}
in agreement with an old classical result \cite{Rob03}.
As $s^\alpha$ has a branch point at $s=0$ in the general case or a pole if
$\alpha$ is a negative integer, we use a deformed integration path $\gamma$
that deviates slightly from $s=0$ in such a way that it remains outside
the integration region. If $\alpha$ is not an integer we choose the negative
real half axis as the branch cut line. This means that
\begin{itemize}
\item $ \delta_\nabla^{(\alpha)} (t) = 0$ for $t<t_0$;
\item for $t= t_n \ge t_0$ the integrand has $n+1$ poles, so
\end{itemize}
\begin{equation*}
\begin{split}
\delta_\nabla^{(\alpha)} (t)
&= -\frac{1}{2\pi i }
\oint_\gamma s^\alpha \prod\limits_{k=1}^{n+1}
\left[1 - s\mu_k(t_0)\right]^{-1} ds\\
&= \frac{(-1)^n}{\prod\limits_{k=1}^{n+1}\mu_k(t_0)}
\oint_\gamma s^\alpha \prod\limits_{k=1}^{n+1}
\left[s - \frac{1}{\mu_k(t_0)}\right]^{-1} ds.
\end{split}
\end{equation*}
We have several different situations between the two extreme cases.
\begin{itemize}
\item \textit{Different graininess values}.
In this case the integrand has $n+1$ simple poles. The corresponding residues
are easily computed:
\begin{equation*}
\delta_\nabla^{(\alpha)} (t) = (-1)^n \sum_{k=1}^{n+1}\mu_k(t_0)^{-\alpha + n}
\prod\limits_{m=1; m\ne k}^{n+1} \left[ \frac{1}{\mu_m(t_0) - \mu_k(t_0)} \right]
\epsilon(t).
\end{equation*}
\item \textit{Constant graininess $\mu_n(t_0) = h$}.
In this case we have a pole with
order $n+1$ at $s=1/h$. The residue theorem allows us to write that
\[
\delta_\nabla^{(\alpha)} (t) = (-1)^n h^{-n-1} \frac{1}{n!}\left.
\frac{d^n s^\alpha}{ds^n}\right|_{s=\frac{1}{h}},
\]
which gives
\begin{equation}
\label{LTZ_61}
\delta_\nabla^{(\alpha)} (t) = h^{-\alpha} \frac{(-\alpha)_n}{n!}\epsilon(nh)
\end{equation}
and leads to the discrete version of the Gr\"unwald--Letnikov
fractional derivative recently introduced in \cite{OrCoTr15},
confirming the coherence of our results. When $h \rightarrow 0$,
we recover the Liouville derivative \cite{Orti11}.
\end{itemize}

In passing, we obtained results useful in defining polynomials on time scales.
Similarly to the ordinary case, we define the ``power function''
of degree $N$ as the $N$th order anti-derivative of the Heaviside unit step
multiplied by $N!$ (cf. \cite{OrCoTr15}). As the Heaviside unit step has
nabla Laplace transform $1/s$, the $N$th order power has nabla Laplace transform
$N!s^{-N-1}$. In fact, with $\alpha = -N$, $N\in \mathbb{Z}_0^+$, we obtain
the ``power function'', agreeing with the time scale and derivative at hand.
For the nabla case, we can use the above results as follows.
\begin{itemize}
\item \textit{Different graininess values}. In this case we obtain
\begin{equation*}
\begin{split}
\left[\epsilon(t)\right]^{\nabla^{-N}} &= \delta_\nabla^{(-N-1)}(t)\\
&= (-1)^n \sum_{k=1}^{n+1}\mu_k(t_0)^{N + n}
\prod\limits_{m=1; m\ne k}^{n+1} \left[ \frac{1}{\mu_m(t_0) - \mu_k(t_0)} \right].
\end{split}
\end{equation*}
\item \textit{Constant graininess $\mu_n(t_0) = h$}.
This case was considered in \cite{OrCoTr15}.
It is easy to obtain from \eqref{LTZ_61}:
\begin{equation*}
\left[\epsilon(nh)\right]^{\nabla^{-N}}
= h^{N} \frac{(N)_n}{n!} = h^{N} \frac{(N+n)!}{N!n!}\epsilon(nh).
\end{equation*}
\end{itemize}

\begin{remark}
The polynomials we obtain here are ``causal'' in the sense
they are nonnull for $t\ge t_0$. For the delta case we
obtain ``anti-causal'' polynomials. These ``power functions''
are univocally defined, contrarily to those obtained
by the approach followed in \cite{BP01}.
\end{remark}


\subsection{Initial and final value theorems}

We now consider the nabla Laplace transform of causal functions
\[
g(t) = f(t)\epsilon(t).
\]
Using the nabla Laplace transform definition, particularized 
to this case, we verify that, in general, it is not possible to obtain an
\emph{Initial Value Theorem}. We only have to remark that the summation goes
from $0$ to $\infty$ and all the terms except one ($\nu_0 \cdot g(0)$) 
depend on $s$ that is inside the region $R_c$. So there is
no value of $s$ that could make zero all terms. For example,
in the constant graininess case, we have
\[
g(0) = \lim\limits_{s \rightarrow \frac{1}{h}} sF(s).
\]
To deduce the \emph{Final Value Theorem} we only have to define
$y(t_n) = h \sum\limits_{k=0}^{n} g(t_k) $. As it is easy to conclude,
if $g(t)$ is a causal function, then $y(t)$ is the convolution of $g(t)$
with the unit step. So the transform of $y(t) $ is $ \frac{G(s)}{s} $.
As it is easy to see,
\begin{equation*}
y(\infty) = \lim\limits_{s \rightarrow 0} G(s) 
= \lim\limits_{s \rightarrow 0} sY(s),
\end{equation*}
which is equal to the usual theorem.


\subsection{Backward compatibility}

Accordingly to what we just said, when the graininess is the constant $h$,
we recover all the results presented in \cite{OrCoTr15}. In particular,
we make two variable transformations: $s = \frac{1 - z^{-1}}{h}$
in the nabla case and $s = \frac{z-1}{h}$ in the delta case.
\begin{itemize}
\item With the above transformations both exponentials degenerate into the
current exponential: $z^n$.
\item Both transforms recover the classic $Z$ transform.
\item Putting $ s= e^{i\omega} $ we obtain the discrete-time Fourier transform.
\end{itemize}
Going further, and taking the limit as $h \rightarrow 0$,
\begin{itemize}
\item The derivatives convert into the classic derivatives. The negative
integer order gives a well-known formulation of the Riemann integral.
\item The exponentials degenerate into one: $e^{st} $.
\item The two nabla Laplace transforms degenerate also into the classic
two-sided Laplace transform.
\item Putting $s= i\omega$, we obtain the Fourier transform.
\end{itemize}


\section{Generalising the transforms}

In Section~\ref{nde} we presented the general forms assumed
by the exponentials: \eqref{LTZ_23} and \eqref{LTZ_25}. With these
we are able to generalise the nabla and delta Laplace transforms.
It is a simple task to rewrite both transforms for a general time scale.
In the nabla case we have
\begin{equation*}
F_\nabla(s)
= \int\limits_{-\infty}^{+\infty} \nu (t) f(t) e_\Delta(t,t_0;-s)\Delta t
= \int\limits_{-\infty}^{+\infty} f(t) e_\Delta(t,t_0;-s)\Delta t
\end{equation*}
and
\begin{equation*}
f(t) = -\frac{1}{2\pi i } \oint_C F_\nabla(s) e_\nabla(\sigma(t),t_0;s) ds;
\end{equation*}
while for the delta case we have, similarly,
\begin{equation*}
F_\Delta(s)
= \int\limits_{-\infty}^{+\infty} \nu (t) f(t) e_\nabla(t,t_0;-s)\nabla t
= \int\limits_{-\infty}^{+\infty} f(t) e_\nabla(t,t_0;-s)\nabla t
\end{equation*}
and
\begin{equation*}
f(t) =  \frac{1}{2\pi i } \oint_C F_\Delta(s) e_\Delta(\rho(t),t_0;s) ds.
\end{equation*}
It is important to remark here that the integration path in the inverse 
transforms becomes the imaginary axis if the graininess is zero at all 
but one value of $t$ in the time scale.


\section{Conclusion}

We introduced a general approach to define exponentials and transforms
on time scales. Starting from the nabla and delta derivatives,
we studied them in parallel and derived general formulae for defining
exponentials as their eigenfunctions. With these exponentials,
we defined two new Laplace transforms and deduced their most important properties.
We obtained existence and unicity, and defined convolution and correlation.
We also considered linear systems and corresponding transfer functions
and impulse response, and defined a general fractional derivative
on time scales from convolution.


\section*{Acknowledgements}

This work was partially supported by National Funds through
the Foundation for Science and Technology of Portugal (FCT),
under project PEst--UID/EEA/00066/2013 (Ortigueira); by CIDMA 
and FCT within project UID\-/MAT/04106/2013 (Torres);
and by project MTM2013-41704-P from the government of Spain (Trujillo).
The authors would like to thank two referees for their valuable
comments and helpful suggestions.




\begin{thebibliography}{xx}

\bibitem{A09}
C. R. Ahrendt, 
{\it The Laplace transform on time scales},
Panamer. Math. J. {\bf 19} (2009), no.~4, 1--36.

\bibitem{AlGr}
A. Aldroubi\ and\ K. Gr\"ochenig, 
{\it Nonuniform sampling and reconstruction in shift-invariant spaces}, 
SIAM Rev. {\bf 43} (2001), no.~4, 585--620. 

\bibitem{AH90}
B. Aulbach\ and\ S. Hilger, 
{\it A unified approach to continuous and discrete dynamics},
in {\it Qualitative theory of differential equations (Szeged, 1988)},
37--56, Colloq. Math. Soc. J\'anos Bolyai, 53 North-Holland, Amsterdam, 1990.

\bibitem{AH90b}
B. Aulbach\ and\ S. Hilger, 
{\it Linear dynamic processes with inhomogeneous time scale},
in {\it Nonlinear dynamics and quantum dynamical systems (Gaussig, 1990)},
9--20, Math. Res., 59 Akademie Verlag, Berlin, 1990.

\bibitem{BaSt}
P. Babu\ and\ P. Stoica, 
{\it Spectral analysis of nonuniformly sampled data -- a review}, 
Digit. Signal Process. {\bf 20} (2010), 359--378.

\bibitem{Bal14}
D. Baleanu, S. Rezapour\ and\ S. Salehi, 
{\it A $k$-dimensional system of fractional finite difference equations}, 
Abstr. Appl. Anal. {\bf 2014} (2014), Art. ID 312578, 8~pp. 

\bibitem{Bas12}
N. R. O. Bastos,
{\it Fractional calculus on time scales},
Ph.D. Thesis, University of Aveiro, 2012.

\bibitem{BMT11}
N. R. O. Bastos, D. Mozyrska\ and\ D. F. M. Torres, 
{\it Fractional derivatives and integrals on time scales 
via the inverse generalized Laplace transform},
Int. J. Math. Comput. {\bf 11} (2011), J11, 1--9.
{\tt arXiv:1012.1555}

\bibitem{MyID:296}
N. Benkhettou, A. M. C. Brito da Cruz\ and\ D. F. M. Torres, 
{\it A fractional calculus on arbitrary time scales: 
fractional differentiation and fractional integration},
Signal Process. {\bf 107} (2015), 230--237.
{\tt arXiv:1405.2813}

\bibitem{BG07}
M. Bohner\ and\ G. Sh. Guseinov,
{\it The convolution on time scales},
Abstr. Appl. Anal. {\bf 2007} (2007),
Art. ID 58373, 24~pp.

\bibitem{BG10}
M. Bohner\ and\ G. Sh. Guseinov,
{\it The Laplace transform on isolated time scales},
Comput. Math. Appl. {\bf 60} (2010), no.~6, 1536--1547.

\bibitem{BG10b}
M. Bohner\ and\ G. Sh. Guseinov,
{\it The $h$-Laplace and $q$-Laplace transforms},
J. Math. Anal. Appl. {\bf 365} (2010), no.~1, 75--92.

\bibitem{BGK11}
M. Bohner, G. Sh. Guseinov\ and\ B. Karpuz,
{\it Properties of the Laplace transform on time scales with arbitrary graininess}, 
Integral Transforms Spec. Funct. {\bf 22} (2011), no.~11, 785--800.

\bibitem{BP01}
M. Bohner\ and\ A. Peterson,
{\it Dynamic equations on time scales},
Birkh\"auser Boston, Boston, MA, 2001.

\bibitem{BP02}
M. Bohner\ and\ A. Peterson,
{\it Laplace transform and $Z$-transform: unification and extension}, 
Methods Appl. Anal. {\bf 9} (2002), no.~1, 151--157.

\bibitem{BP03}
M. Bohner\ and\ A. Peterson,
{\it Advances in dynamic equations on time scales},
Birkh\"auser Boston, Boston, MA, 2003.

\bibitem{duality}
M. C. Caputo,
{\it Time scales: from nabla calculus to delta calculus and vice versa via duality},
Int. J. Difference Equ. {\bf 5} (2010), no.~1, 25--40

\bibitem{C12}
J. L. Cie\'sli\'nski,
{\it New definitions of exponential, hyperbolic and trigonometric functions on time scales},
J. Math. Anal. Appl. {\bf 388} (2012), no.~1, 8--22.

\bibitem{D07}
J. M. Davis, I. A. Gravagne, B. J. Jackson, R. J. Marks, II\ and\ A. A. Ramos,
{\it The Laplace transform on time scales revisited},
J. Math. Anal. Appl. {\bf 332} (2007), no.~2, 1291--1307.

\bibitem{DGM10}
J. M. Davis, I. A. Gravagne\ and\ R. J. Marks, II,
{\it Convergence of unilateral Laplace transforms on time scales}, 
Circuits Systems Signal Process. {\bf 29} (2010), no.~5, 971--997.

\bibitem{DGM10b}
J. M. Davis, I. A. Gravagne\ and\ R. J. Marks, II,
{\it Bilateral Laplace transforms on time scales: convergence, convolution,
and the characterization of stationary stochastic time series},
Circuits Systems Signal Process. {\bf 29} (2010), no.~6, 1141--1165.

\bibitem{FeGr:1}
H. G. Feichtinger\ and\ K. Gr\"ochenig, 
{\it Irregular sampling theorems and series expansions 
of band-limited functions}, 
J. Math. Anal. Appl. {\bf 167} (1992), no.~2, 530--556. 

\bibitem{FeGr:2}
H. G. Feichtinger\ and\ K. Gr\"ochenig, 
Theory and practice of irregular sampling, 
in {\it Wavelets: mathematics and applications}, 
305--363, Stud. Adv. Math, CRC, Boca Raton, FL, 1994.

\bibitem{GlGr}
I. A. Glover and P. M. Grant, 
{\it Digital Communications}, 
Pearson Education Limited, Edinburgh Gate, 
Harlow, Essex, CM20 2JE, England, 3rd edition, 2010.

\bibitem{Gr}
K. Gr\"ochenig, 
{\it A discrete theory of irregular sampling}, 
Linear Algebra Appl. {\bf 193} (1993), 129--150. 

\bibitem{HAQ10}
A. E. Hamza\ and\ M. A. Al-Qubaty,
{\it A remark on the exponential matrix functions on time scales},
Assiut Univ. J. Math. Comput. Sci. {\bf 39} (2010), no.~2, 21--30.

\bibitem{HAQ12}
A. E. Hamza\ and\ M. A. Al-Qubaty,
{\it On the exponential operator functions on time scales},
Adv. Dyn. Syst. Appl. {\bf 7} (2012), no.~1, 57--80.

\bibitem{H90}
S. Hilger,
{\it Analysis on measure chains---a 
unified approach to continuous and discrete calculus},
Results Math. {\bf 18} (1990), no.~1-2, 18--56.

\bibitem{K11}
B. Karpuz,
{\it On uniqueness of the Laplace transform on time scales},
Panamer. Math. J. {\bf 21} (2011), no.~2, 101--110.

\bibitem{MGD08}
R. J. Marks, II, I. A. Gravagne\ and\ J. M. Davis,
{\it A generalized Fourier transform and convolution on time scales},
J. Math. Anal. Appl. {\bf 340} (2008), no.~2, 901--919.

\bibitem{MT09}
D. Mozyrska\ and\ D. F. M. Torres,
{\it A study of diamond-alpha dynamic equations on regular time scales},
Afr. Diaspora J. Math. (N.S.) {\bf 8} (2009), no.~1, 35--47.
{\tt arXiv:0902.1380}

\bibitem{Orti00}
M. D. Ortigueira, 
{\it Introduction to fractional signal processing. 
Part 2: Discrete-time systems}, IEE Proc. on Vision, 
Image and Signal Processing \textbf{147} (2000), no.~1, 71--78.

\bibitem{Orti01}
M. D. Ortigueira,
{\it The comb signal and its Fourier transform},
Signal Processing {\bf 81} (2001), no.~3, 581--592.

\bibitem{Orti11}
M. D. Ortigueira,
{\it Fractional calculus for scientists and engineers},
Lecture Notes in Electrical Engineering, 84, Springer, Dordrecht, 2011.

\bibitem{Orti:14}
M. D. Ortigueira,
{\it On the particular solution of constant coefficient fractional
differential equations},
Appl. Math. Comput. {\bf 245} (2014), 255--260.

\bibitem{OrCoTr15}
M. D. Ortigueira, F.J. Coito\ and\ J. J. Trujillo,
{\it Discrete-time differential systems},
Signal Process. {\bf 107} (2015), 198--217.

\bibitem{Pe}
M. A. Peltola, 
{\it Role of editing of R–R intervals in the analysis of heart rate variability}, 
Front Physiol. {\bf 3} (2012), Art. ID 148, 10~pp. 

\bibitem{PM}
J. G. Proakis\ and\ D. G. Manolakis,
{\it Digital signal processing: Principles, algorithms, and applications},
Prentice Hall, 2007.

\bibitem{Rob03}
M. J. Roberts,
{\it Signals and systems: Analysis using transform methods and Matlab},
McGraw-Hill, 2003.

\bibitem{SR11}
M. R. Segi Rahmat, 
{\it The $(q,h)$-Laplace transform on discrete time scales},
Comput. Math. Appl. {\bf 62} (2011), no.~1, 272--281.

\bibitem{WuBaZe}
G.-C. Wu, D. Baleanu\ and\ S.-D. Zeng, 
Discrete chaos in fractional sine and standard maps, 
Phys. Lett. A {\bf 378} (2014), no.~5-6, 484--487. 

\end{thebibliography}
\end{document}